\documentclass[12pt,notitlepage]{amsart}%
\usepackage[dvips]{graphicx}
\usepackage{pslatex}
\usepackage{pstricks}
\usepackage[myheadings]{fullpage}
\usepackage{mathtools}
\usepackage{tkz-euclide}
\usepackage{subcaption}

\newtheorem{theorem}{Theorem}[section]

\newtheorem{lemma}[theorem]{Lemma}

\newtheorem{proposition}[theorem]{Proposition}

\newcommand{\s}[1]{\sin \frac{#1}{2}}
\newcommand{\sh}[1]{\sinh \frac{#1}{2}}
\newcommand{\f}[1]{\sin \frac{#1}{4}}
\newcommand{\fh}[1]{\sinh \frac{#1}{4}}
\newcommand{\e}[1]{e^\frac{#1}{2}}

\newcommand{\fracsum}{\frac{a}{b}+\frac{b}{c}+\frac{c}{a}}
\newcommand{\fracsym}{\frac{a}{b}+\frac{b}{a}+\frac{a}{c}+\frac{c}{a}+\frac{b}{c}+\frac{c}{b}}
\newcommand{\fracsums}{\frac{\s{a}}{\s{b}}+\frac{\s{b}}{\s{c}}+\frac{\s{c}}{\s{a}}}
\newcommand{\fracsyms}{\frac{\s{a}}{\s{b}}+\frac{\s{b}}{\s{a}}+\frac{\s{a}}{\s{c}}+\frac{\s{c}}{\s{a}}+\frac{\s{b}}{\s{c}}+\frac{\s{c}}{\s{b}}}
\newcommand{\fracsumh}{\frac{\sh{a}}{\sh{b}}+\frac{\sh{b}}{\sh{c}}+\frac{\sh{c}}{\sh{a}}}
\newcommand{\fracsymh}{\frac{\sh{a}}{\sh{b}}+\frac{\sh{b}}{\sh{a}}+\frac{\sh{a}}{\sh{c}}+\frac{\sh{c}}{\sh{a}}+\frac{\sh{b}}{\sh{c}}+\frac{\sh{c}}{\sh{b}}}
    
\begin{document}

\title{Strengthened Euler's Inequality in Spherical and Hyperbolic Geometries}
\author{Ren Guo, Estonia Black, Caleb Smith}
\address{University of Tennessee, Knoxville}
\email{eblack6@vols.utk.edu}
\address{Oregon State University}
\email{smithca6@oregonstate.edu}

\date \today

\thanks{This work was done during the Summer 2016 REU program
in Mathematics at Oregon State University, supported by the National Science Foundation Grant DMS-1359173.}

\begin{abstract}
Euler's inequality is a well known inequality relating the inradius and circumradius of a triangle. In Euclidean geometry, this inequality takes the form $R \geq 2r$ where $R$ is the circumradius and $r$ is the inradius. In spherical geometry, the inequality takes the form $\tan(R) \geq 2\tan(r)$ as proved in \cite{MPV}; similary, we have $\tanh(R) \geq 2\tanh(r)$ for hyperbolic triangles (see \cite{SV} for proof). In Euclidean geometry, this inequality can be strengthened as discussed in \cite{SV}. We prove an analogous version of this strengthened inequality which holds in spherical geometry, as well as an additional strengthening of Euler's inequality which holds in Euclidean geometry and can be generalized into both spherical and hyperbolic geometry.
\end{abstract}

\maketitle
\markboth{Estonia Black, Caleb Smith}{Strengthened Euler's Inequality in Spherical and Hyperbolic Geometries}

\section{Introduction}
In Euclidean geometry, as discussed in \cite{SV}, given a triangle with side lengths $a,b,c$, circumradius $R$, and inradius $r$, we may strengthen Euler's inequality to give:
\begin{theorem}[\cite{SU},\cite{VW}]\label{original}
\begin{subequations}
	\begin{align} 
	\frac{R}{r} &\geq \frac{abc + a^3+b^3+c^3}{2abc} \label{originalA} \\ 
			&\geq \fracsum - 1 \label{originalB} \\ 	
			&\geq \frac{2}{3} \left( \fracsum\right) \label{originalC} \\ 
			&\geq 2 \label{originalD} 
	\end{align} 
\end{subequations} 
\end{theorem}

In this paper, we will discuss how this strengthened formula can be extended to spherical and hyperbolic geometries. We will begin by proving unified results which hold in all three geometries in Section 2. In Section 3, we will prove the strengthened inequality shown above in Euclidean and spherical geometry and discuss how it fails in hyperbolic geometry. Finally, in Section 4 we will discuss an alternative strengthening of the Euler triangle inequality extending our central result in Section 2. 

\section{Unified Results}
In this section, we prove results that unify the three geometries, culminating in Theorem \ref{generalized}, which provides an inequality that holds in each geometry. 
Throughout, we define
\begin{equation*}
	s(x) = \label{s-func}
    \left\{
	\begin{aligned}
						&\frac{x}{2} 	&&\textnormal{in Euclidean geometry} \\
						&\sh{x}			&&\textnormal{in hyperbolic geometry} \\
						&\s{x}			&&\textnormal{in spherical geometry} \\
	\end{aligned}
\right.
\end{equation*}

Our first result will be a lemma, modeled after the central result of \cite{GN}, which will allow us to prove Theorem \ref{spherical}, the spherical analogue of Theorem \ref{original}.


\begin{lemma}\label{unification}
Let $f(a,b,c) \geq 0$ be an inequality which holds for all Euclidean triangles with side lengths $a,b,c$. Then $f(s(a), s(b), s(c)) \geq 0$ for all spherical or hyperbolic triangles with side lengths $a, b, c$, where $s(x)$ is defined as above.

\begin{proof} 
We begin with the hyperbolic case. Consider a hyperbolic triangle with sidelengths $a, b,c$ and vertices $A,B,C$ in the Poincar\'e disk model of the hyperbolic plane. Without loss of generality, we may assume that the triangle is positioned so that its circumcenter coincides with the origin. When we consider the unit disk as a subset of the Euclidean plane, we may consider the Euclidean triangle $T_1$ with vertices $A, B, C$. This triangle will have side lengths $(1-R^2)\sh{a}, (1-R^2)\sh{b}$, and $(1-R^2)\sh{c}$, where $R$ is the Euclidean radius of the circumcircle of $T_1$. This implies that there is another Euclidean triangle $T_2$ similar to $T_1$ with side lengths $\sh{a}, \sh{b}, \sh{c}$. Since the inequality $f(a,b,c) \geq 0$ holds for all Euclidean triangles, we may apply it to $T_2$ to get $f(\sh{a}, \sh{b}, \sh{c}) \geq 0$. This completes the proof for the hyperbolic case.

We now prove the spherical case. Consider a spherical triangle with side lengths $a,b,c$ and vertices $A,B,C$. As in the hyperbolic case, we consider the Euclidean triangle $T_1$ with vertices $A,B,C$. Then the Euclidean distance between $B$ and $C$ will be $2 \s{a}$, since the leg of the Euclidean triangle will be a chord of the great circular arc of length $a$. The other two side lengths can be found in the same fashion. This gives us a Euclidean triangle with side lengths $2\s{a}, 2\s{b}$ and  $2\s{c}$. This implies that there is another Euclidean triangle $T_2$ similar to $T_1$ with side lengths $\s{a}, \s{b}, \s{c}$. Since the inequality $f(a,b,c) \geq 0$ holds for any Euclidean triangle, we may apply it to $T_2$ to get $f(\s{a}, \s{b}, \s{c}) \geq 0$, as required. 

\end{proof}
\end{lemma}

Next we prove a lemma which relates in a unified manner two quantities which appear frequently in inequalities relating $R$ and $r$. This lemma will be instrumental in the proof of Theorem \ref{generalized}.

\begin{lemma}\label{B}
For a triangle in Euclidean, spherical, or hyperbolic geometry, with side-lengths $a,b,c$, let 
\begin{align*}
	\bar{B} \vcentcolon &= \left( s(a)+s(b)-s(c)\right) 	\left( s(a)+s(c)-s(b)\right) \left( s(b)+s(c)-    s(a)\right)\\
	B \vcentcolon &= s(a+b-c)s(a+c-b)s(b+c-a)
\end{align*}
Then
\begin{equation*}
	B-\bar{B}\left\{
	\begin{aligned}
		&=0 		&&\textnormal{in Euclidean geometry} \\
		&\geq 0		&&\textnormal{in hyperbolic geometry} \\
		&\leq 0		&&\textnormal{in spherical geometry} \\
	\end{aligned}
	\right.
\end{equation*}
\begin{proof}
	In Euclidean geometry,
	\begin{align*}
	B-\bar{B}&=\left(\frac{a+b-c}{2}\right)\left(\frac{a+c-b}{2}\right)\left(\frac{b+c-a}{2}\right)-\left(\frac{a}{2}+\frac{b}{2}-\frac{c}{2}\right)\left(\frac{a}{2}+\frac{c}{2}-\frac{b}{2}\right)\left(\frac{b}{2}+\frac{c}{2}-\frac{a}{2}\right)\\
    &=\frac{1}{8}(a+b-c)(a+c-b)(b+c-a)-\frac{1}{8}(a+b-c)(a+c-b)(b+c-a)\\
    &=0.
	\end{align*}
    
    To show that $B\geq\bar{B}$ in hyperbolic geometry, we assume, without loss of generality, that $a\geq b \geq c$.  Then it is sufficient to verify the following two propositions:
    \begin{proposition}
    $\sh{b+c-a}\geq\sh{b}+\sh{c}-\sh{a}$
		\begin{proof}
        Since $2a-b-c\geq b-c$, and $\cosh$ is an even, increasing function, \[\cosh\frac{b+c-2a}{4}=\cosh\frac{2a-b-c}{4}\geq\cosh\frac{b-c}{4},\]
        and therefore
        \[\sh{b+c-a}+\sh{a}=2\fh{b+c}\cosh\frac{b+c-2a}{4}\geq2\fh{b+c}\cosh\frac{b-c}{4}=\sh{b}+\sh{c}\]
		\end{proof}
	\end{proposition}
    \begin{proposition}\label{hypB}
    $\sh{a+b-c}\sh{a+c-b}\geq\left( \sh{a}+\sh{b}-\sh{c}\right) \left( \sh{a}+\sh{c}-\sh{b}\right) $
		\begin{proof}
		$\sh{a+b-c}\sh{a+c-b}\geq\left( \sh{a}+\sh{b}-\sh{c}\right) \left( \sh{a}+\sh{c}-\sh{b}\right)$ if and only if
    	\begin{alignat*}{2}
   		0&\leq&& \sh{a+b-c}\sh{a+c-b}-(\sh{a}+\sh{b}-\sh{c})(\sh{a}+\sh{c}-\sh{b})\\
    	&=&&\frac{1}{4}(\e{a+b-c}-\e{c-a-b})(\e{a+c-b}-\e{b-a-c})-\frac{1}{4}(\e{a}-\e{-a}+\e{b}-\e{-b}-\e{c}+\e{-c})\\& &&\ \ \ (\e{a}-\e{-a}+\e{c}-\e{-c}-\e{b}+\e{-b})\\
    &=&&\frac{1}{4}(e^a+e^{-a}-e^{b-c}-e^{c-b})
    -\frac{1}{4}(2+e^a+e^{-a}-e^b-e^{-b}-e^c-e^{-c}+2(\e{b+c}+\e{-b-c})-2(\e{b-c}+\e{c-b}))\\
    &=&&\frac{1}{4}(e^b+e^{-b}+e^c+e^{-c}+2\e{b-c}+2\e{c-b}-2\e{b+c}-2\e{-b-c}-e^{b-c}-e^{c-b}-2)\\
    &=&&\frac{1}{4}(\e{b}-\e{-b})(\e{b}+\e{2c-b}-2\e{c}-\e{b-2c}-\e{-b}+2\e{-c})\\
    &=&&\frac{1}{2}\sh{b}(\e{c}-\e{-c})(\e{b-c}+\e{c-b}-2)\\
    &=&&2\sh{b}\sh{c}\left(\cosh\left(\frac{b-c}{2}\right)-1\right)\\
    &=&&4\sh{b}\sh{c}\sinh^2\left(\frac{b-c}{4}\right)
    \end{alignat*}
	\end{proof}
    \end{proposition}
    Thus, in hyperbolic geometry, $B=\sh{a+b-c}\sh{a+c-b}\sh{b+c-a}\geq(\sh{a}+\sh{b}-\sh{c})(\sh{a}+\sh{c}-\sh{b})(\sh{b}+\sh{c}-\sh{a})=\bar{B}$
    
    To show that $B\leq\bar{B}$ in spherical geometry, we assume, without loss of generality, that $a\geq b \geq c$.  Then it is sufficient to verify the following two propositions:
    \begin{proposition}
    $\s{b+c-a}\leq\s{b}+\s{c}-\s{a}$
		\begin{proof}
        Since $2a-b-c\geq b-c$, and $\cos$ is a decreasing function on the interval $[0,\pi]$, \[\cos\frac{b+c-2a}{4}=\cos\frac{2a-b-c}{4}\leq\cos\frac{b-c}{4},\]
        and therefore
        \[\s{b+c-a}+\s{a}=2\f{b+c}\cos\frac{b+c-2a}{4}\leq2\f{b+c}\cos\frac{b-c}{4}=\s{b}+\s{c}.\]
		\end{proof}
	\end{proposition}
    \begin{proposition}
    $\s{a+b-c}\s{a+c-b}\leq\left( \s{a}+\s{b}-\s{c}\right) \left( \s{a}+\s{c}-\s{b}\right) $
    \begin{proof}
	This follows directly from the proof of (\ref{hypB}), since \[\sinh ix=i\sin x,\] so 
    \begin{align*}
    &\left( \s{a}+\s{b}-\s{c}\right)\left( \s{a}+\s{c}-\s{b}\right)-\s{a+b-c}\s{a+c-b}\\
    &=i\s{a+b-c}i\s{a+c-b}-\left( i\s{a}+i\s{b}-i\s{c}\right) \left( i\s{a}+i\s{c}-i\s{b}\right) \\
    &=\sh{i(a+b-c)}\sh{i(a+c-b)}-\left(\sh{ia}+\sh{ib}-\sh{ic}\right)\left(\sh{ia}+\sh{ic}-\sh{ib}\right)\\
    &=4\sh{ib}\sh{ic}\sinh^2\frac{i(b-c)}{4}\\
    &=4\s{b}\s{c}\sin^2\frac{b-c}{4}\geq 0
    \end{align*}
    \end{proof}
    \end{proposition}
    Thus, in spherical geometry, 
    \begin{align*}
    B&=\s{a+b-c}\s{a+c-b}\s{b+c-a}\\
    &\leq(\s{a}+\s{b}-\s{c})(\s{a}+\s{c}-\s{b})(\s{b}+\s{c}-\s{a})\\
    &=\bar{B}
    \end{align*}
\end{proof}
\end{lemma}

The equations for inradius and circumradius in Euclidean, hyperbolic, and spherical geometry can be unified as follows: for a triangle with side-lengths $a,b,c$, circumradius $R$, and inradius $r$, 
\begin{equation*}
\frac{2s(a)s(b)s(c)}{\sqrt{s(a+b-c)s(a+c-b)s(b+c-a)s(a+b+c)}}=\left\{
\begin{aligned}
&R   &&\text{in Euclidean geometry}\\
&\tan{R}   &&\text{in spherical geometry}\\
&\tanh{R}   &&\text{in hyperbolic geometry}
\end{aligned}
\right.
\end{equation*}
and
\begin{equation*}
\sqrt{\frac{s(a+b-c)s(a+c-b)s(b+c-a)}{s(a+b+c)}}=\left\{
\begin{aligned}
&r   &&\text{in Euclidean geometry}\\
&\tan{r}   &&\text{in spherical geometry}\\
&\tanh{r}   &&\text{in hyperbolic geometry}
\end{aligned}
\right.
\end{equation*}
The following is a theorem which provides a unified inequality dealing with the inradius and circumradius of triangles in all three of these geometries.
\begin{theorem}\label{generalized}
	A triangle with side-lengths $a,b,c$ has the following property in Euclidean geometry:
	\begin{equation} 
    	\label{gE}\frac{R}{r}\geq \frac{2s(\frac{a+b}{2})s(\frac{a+c}{2})s(\frac{b+c}{2})}{s(a)s(b)s(c)}\geq 2,
    \end{equation}
	while in hyperbolic geometry, 
	\begin{equation} 
    	\label{gH}\frac{\tanh R}{\tanh r}\geq \frac{2s(\frac{a+b}{2})s(\frac{a+c}{2})s(\frac{b+c}{2})}{s(a)s(b)s(c)}\geq 2,
    \end{equation}
	and in spherical geometry,
	\begin{equation} 
    	\label{gS}\frac{\tan R}{\tan r}\geq \frac{2s(\frac{a+b}{2})s(\frac{a+c}{2})s(\frac{b+c}{2})}{s(a)s(b)s(c)}\geq 2.
    \end{equation}
    with equality if and only if $a=b=c$.
	\begin{proof}
   	We begin with the Euclidean case.
    We begin by manipulating the right most inequality.
    \begin{equation*}
    \frac{2s(\frac{a+b}{2})s(\frac{a+c}{2})s(\frac{b+c}{2})}{s(a)s(b)s(c)} = \frac{2 (\frac{a+b}{4})(\frac{a+c}{4})(\frac{b+c}{4})}{\frac{abc}{8}} \geq 2\\ 
    \end{equation*}
    This is equivalent to:
    \begin{align*}
     &\frac{1}{8}(a+b)(b+c)(a+c) \geq abc \\ 
\iff \hspace*{.5in}      &ab^2+a^2b+ac^2+a^2c+bc^2+b^2c + 2abc \geq 8abc \\
\iff  \hspace*{.5in}	 &ab^2+a^2b+ac^2+a^2c+bc^2+b^2c \geq 6abc
    \end{align*}
    
    This final expression follows from the arithmetic-geometric mean inequality since 
    \begin{align*}
    \frac{ab^2+bc^2+ca^2}{3} \geq \sqrt[3]{a^3b^3c^3}
    \iff ab^2+bc^2+ca^2 \geq 3abc
    \end{align*}
    Combining this with the corresponding inequality $a^2b+b^2c+c^2a \geq 3abc$ which can be proved by the same method, we see that indeed 
    \[ ab^2+a^2b+ac^2+a^2c+bc^2+b^2c \geq 6abc \]
    
    This then proves the right inequality of (\ref{gE}).
    
    The left inequality in expression (\ref{gE}) is equivalent to
    \[\frac{2abc}{(a+b-c)(a+c-b)(b+c-a)}\geq \frac{(a+b)(a+c)(b+c)}{4abc},\]
    or
    \[8a^2b^2c^2\geq(a+b-c)(a+c-b)(b+c-a)(a+b)(a+c)(b+c).\]
    
    Now, note that for any circle and any point external to the given circle there are two
        lines tangent to that circle which pass through the point, with the point equidistant 
        from the two points of tangency corresponding to those lines, so we have, for some 
        $x,y,z>0$, that $a=y+z$, $b=x+z$, and $c=x+y$, as depicted in the following diagram:
        
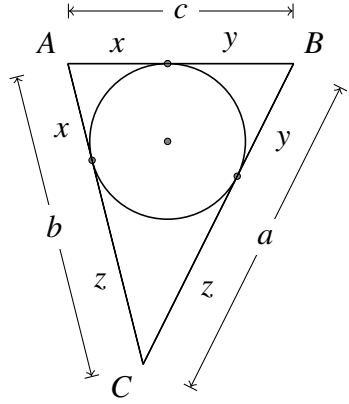
\begin{center}
    	\begin{tikzpicture}

    	\tkzDefPoint(0,2){A}
    	\tkzDefPoint(3,2){B}
    	\tkzDefPoint(1,-2){C}
    	\tkzDrawSegments(A,B B,C C,A)

    	\tkzDefCircle[in](A,B,C)\tkzGetPoint{I}\tkzGetLength{rIN}
    	\tkzDrawPoint(I)
    	\tkzDrawCircle[R](I,\rIN pt)
		\tkzDefPointBy[projection=onto B--C](I) \tkzGetPoint{D} \tkzDrawPoint(D)
    	\tkzDefPointBy[projection=onto A--B](I) \tkzGetPoint{F} \tkzDrawPoint(F)
    	\tkzDefPointBy[projection=onto C--A](I) \tkzGetPoint{E} \tkzDrawPoint(E)
		\draw (A)
		-- (F) node [midway, above]{$x$} -- (B) node [midway, above]{$y$}  -- (A) -- cycle;
		\draw[|<->|] ($(A)!7mm!90:(B)$)--node[fill=white] {$c$} ($(B)!7mm!-90:(A)$);
    	\draw (B)
		-- (D) node [midway, below right]{$y$} -- (C) node [midway, below right]{$z$}  -- (B)
        -- cycle;
		\draw[|<->|] ($(C)!7mm!-90:(B)$)--node[fill=white] {$a$} ($(B)!7mm!90:(C)$);
    	\draw (C)
		-- (E) node [midway, below left]{$z$} -- (A) node [midway, below left]{$x$}  -- (C) -
    	- cycle;
		\draw[|<->|] ($(A)!7mm!-90:(C)$)--node[fill=white] {$b$} ($(C)!7mm!90:(A)$);
    	\tkzLabelPoints[above right](B)
    	\tkzLabelPoints[below left](C)
	    \tkzLabelPoints[above left](A)
		\end{tikzpicture}
        \end{center}
    Then the left inequality of (\ref{gE}) is equivalent to 
    \[8(x+y)^2(x+z)^2(y+z)^2-8xyz(2x+y+z)(x+2y+z)(x+y+2z)\geq0,\]
    or 
    \[
    (x + y + z) (x^3 y^2 + x^2 y^3 - 2 x^2 y^2 z + x^3 z^2 - 
   2 x^2 y z^2 - 2 x y^2 z^2 + y^3 z^2 + x^2 z^3 + y^2 z^3)\geq 0\]
   Which is true, since 
   \begin{align*}
   x^3 y^2 &+ x^2 y^3 - 2 x^2 y^2 z + x^3 z^2 - 2 x^2 y z^2 - 2 x y^2 z^2 + y^3 z^2 + x^2 z^3 + y^2 z^3\\
   &=x^2(y+z)(y-z)^2+y^2(x+z)(x-z)^2+z^2(x+y)(x-y)^2\\
   &\geq 0
   \end{align*}
   for all non-negative $x,y,z$.
    
  	In spherical geometry, the left inequality of (\ref{gS}) is equivalent to
    \[\frac{2\s{a}\s{b}\s{c}}{B}\geq \frac{2\f{a+b}\f{a+c}\f{b+c}}{\s{a}\s{b}\s{c}},\]
    or 
    \[(\s{a}\s{b}\s{c})^2\geq B\f{a+b}\f{a+c}\f{b+c}.\]   
    As discussed in \cite{GN}, we have
    \[\sin^2R=\frac{4(\s{a}\s{b}\s{c})^2}{\bar{B}(\s{a}+\s{b}+\s{c})}\]
    \[\tan^2R=\frac{4(\s{a}\s{b}\s{c})^2}{B\s{a+b+c}}\]
    so
    \begin{align*}
 B\s{a+b+c}+4(\s{a}\s{b}\s{c})^2&=\frac{4(\s{a}\s{b}\s{c})^2}{\tan^2R}+4(\s{a}\s{b}\s{c})^2\\
 	&=\frac{4(\s{a}\s{b}\s{c})^2}{\sin^2R}\\
    &=\bar{B}(\s{a}+\s{b}+\s{c}),
        \end{align*}
     since \[\frac{1}{\sin^2R}=\frac{1}{\tan^2R}+1.\]
     Now, since $B\leq\bar{B}$,    
    \[B\left( \s{a}+\s{b}+\s{c}-\s{a+b+c}\right) \leq4(\s{a}\s{b}\s{c})^2\]
    but since
    \[\s{a}+\s{b}+\s{c}-\s{a+b+c}=4\left( \f{a+b}\f{a+c}\f{b+c}\right)\]
        this is precisely equivalent to the left inequality of (\ref{gS}).
        Now, the right inequality of (\ref{gS}) is equivalent to 
        \[\f{a+b}\f{a+c}\f{b+c}\geq\s{a}\s{b}\s{c},\]
        which is shown on page 636 of \cite{MPV}.
        
    We now consider the hyperbolic case. Similar to the spherical case, the left most inequality in (\ref{gH}) is equivalent to 
    \begin{align*}
    \frac{2\sh{a}\sh{b}\sh{c}}{B}\geq \frac{2\fh{a+b}\fh{a+c}\fh{b+c}}{\sh{a}\sh{b}\sh{c}} \\
    \iff (\sh{a}\sh{b}\sh{c})^2 \geq B \fh{a+b} \fh{a+c} \fh{b+c}
	\end{align*}
    Just as in the spherical case, we have expressions for 
    \begin{align*}
    \sinh^2 R &= \frac{4(\sh{a}\sh{b}\sh{c})^2}{\bar{B}(\sh{a}+\sh{b}+\sh{c})} \\
    \tanh^2 R &= \frac{4(\sh{a}\sh{b}\sh{c})^2}{B \sh{a+b+c}}
    \end{align*}
    as discussed in \cite{GN}. Using the formula $\frac{1}{\sinh^2x} = \frac{1}{\tanh^2 x} -1$, we get
    \begin{align*}
    \bar{B}\left(\sh{a}+\sh{b}+\sh{c}\right) = B \sh{a+b+c} - 4\left(\sh{a}\sh{b}\sh{c}\right)^2
    \end{align*}
    Using the fact that $B \geq \bar{B}$ and \[\sh{a+b+c}-\sh{a}-\sh{b}-\sh{c} = 4\fh{a+b}\fh{a+c}\fh{b+c},\] we are left with
    \begin{equation*}
    B\fh{a+b}\fh{a+c}\fh{b+c} \leq \left(\sh{a}\sh{b}\sh{c}\right)^2
    \end{equation*}
    which, as was shown above, is equivalent to the leftmost inequality of expression (\ref{gH}).
   \end{proof}
\end{theorem}

\section{Original Strengthening}

\subsection{Euclidean}

We will here prove Theorem \ref{original} for the reader's convenience.
\begin{proof}
The proof of (\ref{original}) will be divided into proofs of its four parts:
\begin{itemize}
\item[(\ref{originalD})] $\frac{2}{3}(\fracsum)\geq2$:

		(\ref{originalD})$\iff\fracsum\geq3$. Now, by the inequality of arithmetic and geometric means, 
        \[\frac{1}{3}\left( \fracsum\right) \geq\sqrt[3]{\frac{a}{b}\frac{b}{c}\frac{c}{a}}=1,\]
        so 
        $\fracsum\geq3$

\item[(\ref{originalC})] $\fracsum-1\geq\frac{2}{3}(\fracsum)$:

    	This follows directly from the inequality proven above: 
        $\frac{2}{3}(\fracsum)\geq2\iff\frac{1}{3}(\fracsum)\geq1\iff\fracsum-1\geq\fracsum-\frac{1}{3}(\fracsum)=\frac{2}{3}(\fracsum)$.
\item[(\ref{originalB})] $\frac{abc+a^3+b^3+c^3}{2abc}\geq\fracsum-1$:

		$\frac{abc+a^3+b^3+c^3}{2abc}\geq\fracsum-1=\frac{a^2c+b^2a+c^2b-abc}{abc}$, or  $abc+a^3+b^3+c^3\geq 2(a^2c+b^2a+c^2b-abc)$ if and only if
        \begin{align*}
        0&\leq a^3+b^3+c^3+3abc-2a^2c-2b^2a-2c^2b\\
        &=(a^3+ac^2-2a^2c)+(b^3+a^2b-2ab^2)+(c^3+b^2c-2bc^2)-ac^2-a^2b-b^2c+3abc\\ 
        &=a(a-c)^2+b(b-a)^2+c(c-b)^2+ac(b-c)+ab(c-a)+bc(a-b).  
        \end{align*}
        Without loss of generality, suppose $a\geq b\geq c$. Now
		\begin{align*}
  		0&\leq (a-c)^2(a-b)\\
        &=(a-c)(a(a-c)+bc-ab)\\
        &=(a-c)(a^2-ac+bc-ab)\\
        &=a(a-c)^2+bc(b-c)+ab(c-a)+bc(a-b)\\
        &\leq a(a-c)^2+ac(b-c)+ab(c-a)+bc(a-b)\\
        &\leq a(a-c)^2+b(b-a)^2+c(c-b)^2+ac(b-c)+ab(c-a)+bc(a-b)
        \end{align*}
\item[(\ref{originalA})] $\frac{R}{r}\geq\frac{abc+a^3+b^3+c^3}{2abc}$:
       
       Recall that we have $a=y+z$, $b=x+z$, and $c=x+y$ for some $x,y,z>0$.
        Now, \[\frac{R}{r}=\frac{abc}{4(\frac{a+b-c}{2})(\frac{a+c-b}{2})(\frac{b+c-a}{2})}=\frac{2abc}{(a+b-c)(a+c-b)(b+c-a)},\] 
        so
        \[\frac{R}{r}\geq\frac{abc+a^3+b^3+c^3}{2abc}\]
        if and only if
		 \[(y+z)^2(x+z)^2(x+y)^2\geq 2xyz((y+z)(x+z)(x+y)+(y+z)^3+(x+z)^3+(x+y)^3),\]
         which is equivalent to        
        \[2(x^2y^2(x-z)(y-z)+y^2z^2(y-x)(z-x)+x^2z^2(x-y)(z-y))+x^4(y-z)^2+y^4(x-z)^2+z^4(y-z)^2\geq0.\]
        So it will suffice to show 
        \[x^2y^2(x-z)(y-z)+y^2z^2(y-x)(z-x)+x^2z^2(x-y)(z-y)\geq0\]
       Suppose, without loss of generality, that $x\geq y\geq z$. Then $y^2\geq z^2$ and $(x-z)\geq(x-y)$, so 
       \begin{align*}
       0&\leq y^2z^2(y-x)(z-x)\\
       &=y^2z^2(y-x)(z-x)+z^2x^2(y-z)(x-y)+z^2x^2(z-y)(x-y)\\
       &\leq x^2y^2(x-z)(y-z)+y^2z^2(y-x)(z-x)+x^2z^2(x-y)(z-y)
        \end{align*}  
\end{itemize}
\end{proof}

\subsection{Spherical}

\begin{theorem}\label{spherical}
	Let $a,b,$ and $c$ be the side-lengths of a spherical 
    triangle with circumradius $R$ and inradius $r$. Then:
    \begin{subequations}
    	\begin{align}
		\frac{\tan R}{\tan r}&\geq
        \frac{\s{a}\s{b}\s{c}+\sin^3 \frac{a}{2}+\sin^3 \frac{b}{2}+\sin^3{\frac{a}{2}}}{2\s{a}\s{b}\s{c}} 	\label{sphericalA} \\
        &\geq\fracsums-1 	\label{sphericalB} \\
        &\geq\frac{2}{3} \left(\fracsums\right) 	\label{sphericalC} \\
        &\geq2\label{sphericalD}
    	\end{align}
    \end{subequations}
	\begin{proof}
	(\ref{sphericalB}), (\ref{sphericalC}), and (\ref{sphericalD}) follow from (\ref{originalB}), (\ref{originalC}), and (\ref{originalD}), respectively, as an application of Lemma \ref{unification}, so it remains only to show (\ref{sphericalA}).
    
    Now, note that
    \begin{align*}
    \tan R&=\frac{2\s{a}\s{b}\s{c}}{\sqrt{\s{a+b+c}\s{a+b-c}\s{a+c-b}\s{b+c-a}}}\\
	\tan r&=\sqrt{\frac{\s{a+b-c}\s{a+c-b}\s{b+c-a}}{\s{a+b+c}}}
	\end{align*}
	so
	\[\frac{\tan R}{\tan r}=\frac{2\s{a}\s{b}\s{c}}{\s{a+b-c}\s{a+c-b}\s{b+c-a}}=\frac{2\s{a}\s{b}\s{c}}{B},\]
	and since $B\leq\Bar{B}$, we have
    \[\frac{\tan R}{\tan r}\geq\frac{2\s{a}\s{b}\s{c}}{\bar{B}}\geq \frac{\s{a}\s{b}\s{c}+\sin^3 \frac{a}{2}+\sin^3 \frac{b}{2}+\sin^3{\frac{a}{2}}}{2\s{a}\s{b}\s{c}},\]
    by Lemma \ref{unification}.
	\end{proof}
\end{theorem}

\subsection{Hyperbolic} \label{section-hyperbolic-orig}
Although we have 
\begin{align*}
	\frac{\sh{a}\sh{b}\sh{c}+\sinh^3\frac{a}{2}+\sinh^3\frac{b}{2}+\sinh^3\frac{c}{2}}{2\sh{a}\sh{b}\sh{c}}&\geq \fracsumh-1\\
    &\geq\frac{2}{3}\left( \fracsumh\right)\\
    &\geq 2
\end{align*}
for triangles in hyperbolic geometry with side-lengths $a,b,c,$ circumradius $R$, and inradius $r$ as an application of Lemma \ref{unification}, there is no direct generalization of Theorem \ref{original} into hyperbolic geometry. 
Take, for example, a triangle in hyperbolic geometry with side-lengths $a=b=2$ and $c=0.4$.  Then
\begin{align*}
	0&\geq -0.00923\\
    &\approx\frac{\tanh R}{\tanh r}-\frac{2}{3}\left( \fracsumh\right)\\
    &\geq\frac{\tanh R}{\tanh r}-\left( \fracsumh-1\right)\\
    &\geq\frac{\tanh R}{\tanh r}-\frac{\sh{a}\sh{b}\sh{c}+\sinh^3\frac{a}{2}+\sinh^3\frac{b}{2}+\sinh^3\frac{c}{2}}{2\sh{a}\sh{b}\sh{c}}
\end{align*}

In fact, for triangles in hyperbolic geometry, $\frac{\tanh R}{\tanh r}$ is not comparable with a generalized form of any of the terms in Theorem \ref{original}, apart from that given by the original Euler's inequality. To see this, consider a triangle with edge-lengths $a=b=2$ and $c=0.5$. Then we have
\begin{align*}
	0&\leq0.037107\\
    &\approx \frac{\tanh R}{\tanh r}-\frac{\sh{a}\sh{b}\sh{c}+\sinh^3\frac{a}{2}+\sinh^3\frac{b}{2}+\sinh^3\frac{c}{2}}{2\sh{a}\sh{b}\sh{c}}\\
    &\leq\frac{\tanh R}{\tanh r}-\left( \fracsumh-1\right)\\
    &\leq\frac{\tanh R}{\tanh r}-\frac{2}{3}\left( \fracsumh\right)
\end{align*}
\section{Symmetric Strengthening}

\subsection{Euclidean}

\begin{theorem}\label{symmetric}
A Euclidean triangle with edge-lengths $a,b,c,$ circumradius $R$, and inradius $r$ has the following property:
	\begin{subequations}
    \begin{align}
		\frac{R}{r}&\geq\frac{abc+a^3+b^3+c^3}{2abc}\label{symmetricA}\\
    	&\geq \frac{1}{2}\left( \fracsym\right)\label{symmetricB} -1\\
        &\geq \frac{1}{3}\left( \fracsym\right)\label{symmetricC}\\
        &\geq\frac{(a+b)(a+c)(b+c)}{4abc}\label{symmetricD}\\
        &\geq 2\label{symmetricE}
    \end{align}
	\end{subequations}
    \begin{proof}
	(\ref{symmetricA}) is equivalent to (\ref{originalA}). (\ref{symmetricB}) follows from (\ref{originalB}), since
    \begin{align*}
	\frac{abc+a^3+b^3+c^3}{2abc}&\geq \fracsum-1\\
    \frac{abc+a^3+b^3+c^3}{2abc}&\geq \frac{b}{a}+\frac{a}{c}+\frac{c}{b}-1
	\end{align*}
    so
    \[\frac{abc+a^3+b^3+c^3}{2abc}\geq \frac{1}{2}\left( \fracsym -2\right).\]
    Similarly, (\ref{symmetricC}) follows from (\ref{originalC}).
    
    (\ref{symmetricD}) and (\ref{symmetricE}) are both equivalent to 
    \[a^2b+b^2a+a^2c+c^2a+b^2c+c^2b\geq6abc\]
    which is true, since by the inequality of arithmetic and geometric means
    \begin{align*}
    	ab^2+ac^2&\geq2abc\\
        a^2b+bc^2&\geq2abc\\
        a^2c+b^2c&\geq2abc.
    \end{align*}
    \end{proof}
\end{theorem}

\subsection{Spherical}\label{(sphericalsym)}
We have the following theorem for triangles in spherical geometry:
\begin{theorem}
	A triangle in spherical geometry with edge-lengths $a,b,c,$ circumradius $R$, and inradius $r$ has the following property:
	\begin{subequations}
	\begin{align}
    \frac{\tan R}{\tan r}&\geq \frac{\s{a}\s{b}\s{c}+\sin^3\frac{a}{2}+\sin^3\frac{b}{2}+\sin^3\frac{c}{2}}{2\s{a}\s{b}\s{c}}\label{symmetricsA}\\
    &\geq \frac{1}{2}\left( \fracsyms\right) -1\label{symmetricsB}\\
    &\geq \frac{1}{3}\left( \fracsyms\right)\label{symmetricsC}\\
    &\geq 2\label{symmetricsD}
    \end{align}
	\end{subequations}
    \begin{proof}
    	(\ref{symmetricsA}) is equivalent to (\ref{sphericalA}). (\ref{symmetricsB}) and  (\ref{symmetricsC}) follow from (\ref{symmetricB}) and (\ref{symmetricC}), respectively, as an application of Lemma \ref{unification}.
        
        (\ref{symmetricsD}) follows from (\ref{sphericalD}), since 
        \begin{align*}
        	\frac{2}{3}\left( \fracsums\right) &\geq2\\
            \frac{2}{3}\left( \frac{\s{b}}{\s{a}}+\frac{\s{a}}{\s{c}}+\frac{\s{c}}{\s{b}}\right) &\geq2
        \end{align*}
        so
        \[\frac{2}{3}\left( \fracsyms\right) \geq 4\]
        \end{proof}
\end{theorem}

We cannot, however, include an inequality analogous to (\ref{symmetricD}) in spherical geometry. Take, for example, a triangle in spherical geometry with edge-lengths $a=b=3$ and $c=1.5$. Then,
\begin{align*}
	0&\leq 0.19775\\
	&\approx \frac{2\sin\frac{a+b}{4}\sin\frac{a+c}{4}\sin\frac{b+c}{4}}{\s{a}\s{b}\s{c}}-\frac{\s{a}\s{b}\s{c}+\sin^3\frac{a}{2}+\sin^3\frac{b}{2}+\sin^3\frac{c}{2}}{2\s{a}\s{b}\s{c}}\\
	&\leq\frac{2\sin\frac{a+b}{4}\sin\frac{a+c}{4}\sin\frac{b+c}{4}}{\s{a}\s{b}\s{c}}-\left(\frac{1}{2}(\fracsyms)-1\right)\\
    &\leq \frac{2\sin\frac{a+b}{4}\sin\frac{a+c}{4}\sin\frac{b+c}{4}}{\s{a}\s{b}\s{c}}-\frac{1}{3}\left( \fracsyms\right).
\end{align*}
In fact, $\frac{2\sin\frac{a+b}{4}\sin\frac{a+c}{4}\sin\frac{b+c}{4}}{\s{a}\s{b}\s{c}}$ is not comparable with any of the quantities expressed in Theorem \ref{(sphericalsym)}, except for those given in Theorem \ref{generalized}. To see this, consider a triangle with side-lengths $a=b=0.75$, $c=1$, which has
\begin{align*}
	0&\leq0.00418\\
    &\approx\frac{1}{3}\left( \fracsyms\right)-\frac{2\f{a+b}\f{a+c}\f{b+c}}{\s{a}\s{b}\s{c}}.\\
    &\leq\frac{1}{2}\left( \fracsyms\right)-1-\frac{2\f{a+b}\f{a+c}\f{b+c}}{\s{a}\s{b}\s{c}}\\
    &\leq\frac{\s{a}\s{b}\s{c}+\sin^3\frac{a}{2}+\sin^3\frac{b}{2}+\sin^3\frac{c}{2}}{2\s{a}\s{b}\s{c}}-\frac{2\f{a+b}\f{a+c}\f{b+c}}{\s{a}\s{b}\s{c}}.
\end{align*}
\subsection{Hyperbolic}
Although we have
\begin{align*}
	&\frac{\sh{a}\sh{b}\sh{c}+\sinh^3\frac{a}{2}+\sinh^3\frac{b}{2}+\sinh^3\frac{c}{2}}{2\sh{a}\sh{b}\sh{c}}\\
    &\geq \frac{1}{2}\left( \fracsymh\right) -1\\
	&\geq \frac{1}{3}\left( \fracsymh\right)\\
	&\geq 2
\end{align*}
for triangles of edge-lengths $a,b,c$ in hyperbolic geometry as an application of Lemma \ref{unification} to Theorem \ref{symmetric}, there is no theorem in hyperbolic geometry which is directly analogous to Theorem \ref{symmetric}.
To show this, consider a triangle with side-lengths $a = b = 2.5$ and $c = 2$, then we have
\begin{align*}
 0&\geq-0.0457201\\ 
 &\approx \frac{\tanh{R}}{\tanh{r}} - \frac{1}{3} \left(\fracsymh \right) \\
 &\geq \frac{\tanh{R}}{\tanh{r}}-\left( \frac{1}{2}(\fracsymh)-1\right)\\
 &\geq \frac{\tanh{R}}{\tanh{r}}-\frac{\sh{a}\sh{b}\sh{c}+\sinh^3\frac{a}{2}+\sinh^3\frac{b}{2}+\sinh^3\frac{c}{2}}{2\sh{a}\sh{b}\sh{c}}.
\end{align*}
As we saw with the failure of Theorem \ref{original} to generalize into hyperbolic geometry, $\frac{\tanh R}{\tanh r}$ is not in fact comparable with any of these quantities, since, for example, a triangle with edge-lengths $a=b=1$, $c=1.5$ has
\begin{align*}
	0&\leq0.23557\\
	&\approx\frac{\tanh R}{\tanh r}-\frac{\sh{a}\sh{b}\sh{c}+\sinh^3\frac{a}{2}+\sinh^3\frac{b}{2}+\sinh^3\frac{c}{2}}{2\sh{a}\sh{b}\sh{c}}\\
    &\leq\frac{\tanh{R}}{\tanh{r}}-\left( \frac{1}{2}(\fracsymh)-1\right)\\
    &\leq\frac{\tanh{R}}{\tanh{r}} - \frac{1}{3} \left(\fracsymh \right).
\end{align*}

\newpage
\bibliographystyle{plain}
\bibliography{bib.bib}

\end{document}